\documentclass[12pt]{amsart}
\usepackage{amsmath,amssymb,amscd,amsthm,amsxtra}
\usepackage{latexsym}
\usepackage{graphicx}

\usepackage[colorlinks=true, pdfstartview=FitV, linkcolor=blue, citecolor=blue, urlcolor=blue]{hyperref}
\usepackage{color}

\usepackage{graphicx,epstopdf,verbatim}

\allowdisplaybreaks

\numberwithin{equation}{section}

\newtheorem{thm}{THEOREM}[section]
\newtheorem{prop}[thm]{PROPOSITION}
\newtheorem{lemma}[thm]{LEMMA}
\newtheorem{coroll}[thm]{COROLLARY}

\newcommand{\strt}[1]{\rule{0pt}{#1}}

\newcommand{\norm}[1]{\mbox{$\left\| #1 \right\|$}}
\newcommand{\absval}[1]{\mbox{$|#1|$}}

\newcommand{\R}{ \mathbb{R} }
\newcommand{\RN }{\mathbb R^n}

\let \e=\varepsilon
\let \d=\delta

\begin{document}

\title[ Quadratic $A_1$ bounds for commutators]{Quadratic $A_1$ bounds for commutators of singular
integrals with BMO functions}

\author{Carmen Ortiz-Caraballo}

\begin{abstract}

For any Calder\'{o}n--Zygmund operator $T$ and any $BMO$
function $b$ we prove the following quadratic estimate
$$
\norm{[b,T]}_{L^p(w)}\le c \norm{b}_{BMO} (pp')^2[w]_{A_1}^2, \qquad
1<p<\infty, \,w\in A_1,
$$
with constant $c=c(n,T)$ being the estimate optimal
on $p$ and the exponent of the weight constant.  As an endpoint estimate we
prove
$$
w(\{x\in \mathbb R^n:|[b,T]f(x)|>\lambda\}\leq c \,\Phi([w]_{A_1})^2
\int_{\mathbb R^n} \Phi \left({{|f(x)|}\over {\lambda}}
\right)w(x)\,dx,
$$
where $\Phi(t)=t(1+\log^{+}t)$ and constant $c=c(n,T,\|b\|_{{BMO}})$.

\end{abstract}

\address{
Departamento de Matem\'aticas
\\ Escuela Polit\'ecnica, Universidad de Extremadura\\ Avda. Universidad, s/n, 10003 C\'aceres, Spain} \email{carortiz@unex.es}



\maketitle

%%%%%%%%%%%%%%%%%%%%%%%%%%%%%%%%%%%%%%%%%%%%%%%%%%%%%%%%%%%

\section{Introduction and main results}

The purpose of this paper is to study $A_1$ weight estimates for
commutators of singular integrals with $BMO$ functions. Before
stating the main results we will describe briefly some of the
motivations and earlier developments.

In 1971, C. Fefferman and E.M. Stein \cite{FS} established the
following extension of the classical weak-type $(1,1)$ property of
the Hardy-Littlewood maximal operator $M$:
\begin{equation}\label{fs}
\norm{M}_{L^{1,\infty}(w)}  \le c\, [w]_{A_1},  \qquad w\in A_1
\end{equation}
Recall that $w$ is an $A_1$ weight if there is a finite constant $c$
such that
$$Mw\le c\,w \quad \mbox{a.e.}$$
and we denote by $[w]_{A_1}$ the smallest of these $c$.
In fact they proved something better:
\begin{equation}\label{fs2}
\norm{Mf}_{L^{1,\infty}(w)} \le c\, \int_{\R^n}|f|\,Mwdx, \qquad
w\geq 0.
\end{equation}
which was used to derive vector-valued extensions of the classical
estimates for $M$.

It was conjectured by B. Muckenhoupt and R. Wheeden in the seventies
that the analogue of \eqref{fs} holds for $T$, namely,
\begin{equation}\label{conjMW}
\norm{Tf}_{L^{1,\infty}(w)} \le c\, \int_{\R^n}|f|\,Mwdx, \qquad
w\geq 0,
\end{equation}
but has been disproved for T the Hilbert transform in a recent paper by Mar\'{\i}a del
Carmen Reguera and Christoph Thiele \cite{RgT}.

The best result in this sense can be found in \cite{Pe3} where $M$ is
replaced by $M_{L(\log L)^{\varepsilon}} $ for any $\varepsilon> 0$
with new constant $c_{\varepsilon}$ which blows up when $\varepsilon
\to 0$. Furthermore, the following weaker
variant of \eqref{conjMW}:
\begin{equation}\label{w-conj}
\norm{H}_{L^{1,\infty}(w)}  \le c\, [w]_{A_1},     \qquad  w\in A_1,
\end{equation}
seems to be false (see  \cite{NRVV}).

Of course, the corresponding weak $L^p$ type estimates have also a
lot of interest. If we consider again the Hardy-Littlewood maximal
function $M$, an application of the Marcinkiewicz Interpolating
Theorem gives
\begin{equation}
\label{weakannonimous} \|M\|_{L^p(w)} \leq c_{n}p'\,[w]_{A_1}^{1/p}
\qquad w\in A_1
\end{equation}
However, similar estimates for the Hilbert transform are more
difficult to prove as shown by R. Fefferman and J. Pipher in
\cite{FPi}. In this paper the authors established the following
result
\begin{equation}\label{a1}
\|H\|_{\strt{1.7ex} L^p(w)}\le c_{n,p}\,[w]_{A_1} \qquad p\geq 2
\end{equation}
being the exponent of $[w]_{A_1}$ best possible. This estimate was
improved in several directions by A. Lerner, S. Ombrosy and C.
P\'{e}rez, in \cite{LOPe1} and \cite{LOPe3} (see also \cite{LOPe2} for a dual problem). Indeed, this result was extended
to any $1<p<\infty$ and to  any Calder\'on-Zygmund operator. Moreover
an endpoint estimate close to the conjecture
\eqref{w-conj} was also obtained. To be more precise, they proved the following
results.

\begin{thm}\cite{LOPe3}\label{lineargrowth}
Let $T$ be a Calder\'on-Zygmund operator and let $1<p<\infty$. Then
\begin{equation} \label{sharp1}
\|T\|_{L^p(w)}\le c\,pp'\, [w]_{A_1}
\end{equation}
where $c=c_{n,T}$. Furthermore this result is optimal.
\end{thm}

We estate now the corresponding result related  to conjecture \eqref{w-conj}.

\begin{thm}\cite{LOPe3}\label{logarithmicgrowth}

Let $T$ be a Calder\'on-Zygmund operator. Then
\begin{equation}\label{improve}
\|T\|_{L^{1,\infty}(w)}\le c\,\Phi([w]_{A_1}),
\end{equation}
where $c=c_{n,T}$ and where $\Phi(t)=t(1+\log^{+}t)$.
\end{thm}

These results were motivated by important works by S. Petermichl and
A. Volberg \cite{PetV} for the Ahlfors-Beurling transform and by S.
Petermichl \cite{Pet1,Pet2} for the Hilbert transform and the Riesz
Transforms. In these papers it has been shown that if $T$ is any of
these operators, then
\begin{equation}\label{pet}
\|T\|_{\strt{1.7ex}L^p(w)}\le
c_{p,n}\,[w]_{A_p}^{\max\{1,\frac{1}{p-1}\}}\qquad 1<p<\infty,
\end{equation}
where the exponent $\max\{1,\frac{1}{p-1}\}$ is best possible. Note
that $A_1\subset A_p$, and $[w]_{A_p}\le [w]_{A_1}$ (See Section
\ref{preliminares} for the definition of $[w]_{A_p}$). Therefore
\eqref{pet} clearly gives the right exponent in \eqref{a1} when
$p\ge 2$. However, (\ref{pet}) cannot be used in order to get the
sharp exponent in the range $1<p<2$, becoming the exponent worst
when $p$ gets close to 1. We note also that the proofs in
\cite{Pet1,Pet2,PetV} are based on the Bellman function techniques for $p=2$. The case $p\neq2$ follows by the sharp version of the extrapolation theorem of Rubio de Francia as can be found in  \cite{DGrPPet} , and
it is not clear whether they can be extended to the wider class of
Calder\'on-Zygmund operators  as is done in
\cite{LOPe1,LOPe3} in the case of $A_1$. We remit the reader to
\cite{Pe4} for a survey on this topic and to the papers \cite{CrMPe1}  and \cite{CrMPe2} for a recent and different approach to these problems.

The sharp bound \eqref{pet}  for any Calder\'on-Zygmund operator $T$ has been proved recently in \cite{Hyt} by T. Hyt\"onen.
As before, it is enough to consider  the case  $p=2$ by the sharp extrapolation theorem. Hyt\"onen's proof is based on  approximating $T$ by generalized dyadic Haar shift operators with good bounds combined with the key fact that to prove \eqref{pet} it is enough to prove the weak type $(2,2)$ estimate with the same linear bound as proved in \cite{PeTV}. A direct proof avoiding this weak $(2,2)$ reduction can be found in \cite{HytPeTV}.  A bit earlier, in \cite{L2}, the sharp $L^p(w)$ bound for $T$ was obtained for values of $p$ outside the interval $(3/2,3)$ and the proof is based on the corresponding estimates for the intrinsic square function. Even more recently, the bound \eqref{pet}  for any Calder\'on-Zygmund has been
further improved by T. Hyt\"onen and C. P\'erez in \cite{HytPe}, where a portion of the $A_p$ constant of $w$ is replaced by the weaker $A_{\infty}$ constant.

The main purpose of this paper is to prove similar estimates to
Theorems \ref{lineargrowth} and \ref{logarithmicgrowth} for
commutators of singular integral operators $T$ with $BMO$ functions
$b$. These operators were introduced by Coifman, Rochberg and Weiss
in \cite{CRoW} and formally they are defined by
\begin{equation}
[b,T]f(x)=b(x) T(f)(x) - T(b\,f)(x) =\int_{\mathbb R^n}
(b(x)-b(y))K(x,y)f(y)\,dy,
\end{equation}
where $K$ is a kernel satisfying the standard Calder\'on-Zygmund
estimates (see Section \ref{preliminares} for the precise
definition). Although the original interest in the study of such
operators was related to generalizations of the classical
factorization theorem for Hardy spaces many other applications  have
been found and in particular in partial differential equations.

The main result from \cite{CRoW} states that $[b, T]$ is a bounded
operator on $L^{p}( \mathbb R^{n} )$, $1<p<\infty$, when $b$ is a
$BMO$ function. In fact, the  $BMO$ condition of $b$ is also a
necessary condition for the $L^{p}$-boundedness of the commutator
when $T$ is the Hilbert transform. These operators often behave as
Calder\'on-Zygmund operators but there are some differences. For
instance, an interesting fact is that, unlike what it is done with
singular integral operators, the proof of the $L^{p}$-boundedness of
the commutator does not rely on a weak type $(1,1)$ inequality.  In
fact, simple examples show that in general $[b,T]$ fails to be of
weak type $(1,1)$ when $b \in BMO$. This was observed by P\'erez in
\cite{Pe1} where it is also shown that there is an appropriate
weak-$L(\log L)$ type estimate replacement (see below). Also it is
shown by the same author in \cite{Pe2} that $M$ is not the right
operator controlling $[b,T]$ but $M^2=M\circ M$. These results
amount to say that they behave differently from the
Calder\'on-Zygmund operators.

In the present  paper we pursue this point of view by showing that
commutators have an extra ``bad" behavior from the point of view of
$A_1$ weights when trying to derive theorems such as Theorems
\ref{lineargrowth} or \ref{logarithmicgrowth}. Related to the first
theorem we have the following result.

\begin{thm}\label{fuerte}
Let $T$ be a Calder\'{o}n--Zygmund operator and let $b$ be in $BMO$.
Also let $1<p,r<\infty$. Then there exists a constant $c=c_{n,T}$
such that for any weight $w$, we claim that the following inequality
holds
\begin{equation}\label{acotadospesos}
\norm{[b,T]f}_{L^{p}(w)}\le c\,  \norm{b}_{BMO}\,
{(pp')}^2\,(r')^{1+\frac{1}{p'}}\,
   \norm{f}_{L^{p}(M_{r}w)}.
\end{equation}
In particular if $w \in A_1$, we have
\begin{equation}
\norm{[b,T]}_{L^p(w)}\le c\,\norm{b}_{BMO} (pp')^2[w]_{A_1}^2.
\end{equation}
Furthermore this result is sharp in $p$ and in the exponent of  $[w]_{A_1}$.
\end{thm}

It should be mentioned that D. Chung proved in his dissertation (2010) that the commutator $[b,H]$, where  $H$ is the Hilbert transform and $b\in BMO$, obeys a quadratic bound in $L^2(w)$ with respect to the $A_2$ constant of the weight \cite{Ch}. His proof is based on  dyadic methods combined with Bellman functions techniques. On the other hand there is a new proof following an idea from \cite{CRoW} by Chung-Pereyra-P\'erez \cite{ChPP}. This result is more general  and states that if a linear operator $T$ which obeys a linear bound in $L^2(w)$ with respect to the
$A_2$ constant, then its corresponding commutator obeys a quadratic bound. In light of Hyt\"onen's
result this implies that all commutators of Calder\'on-Zygmund singular integral operators and
BMO functions obey a quadratic bound in $L^2(w)$ which can then be extrapolated to $L^p(w)$. These
results have been generalized by Cruz-Uribe and Moen \cite{CrMo} to the two-weight setting, to fractional
integrals, and to vector-valued extensions.

The second main result of this paper is the following endpoint
version of Theorem \ref{fuerte}.

\begin{thm}\label{debil}
Let $T$ and $b$ as above. Then there exists a constant $c=c_{n,T,\norm{b}_{BMO}}$
such that for any weight $w \in A_1$ and $f \in L_c^{\infty}(\mathbb
R^n)$
\begin{equation}
w(\{x\in \mathbb R^n:|[b,T]f(x)|>\lambda\}\leq c\,
\Phi([w]_{A_1})^2 \int_{\mathbb R^n} \Phi
\left({{|f(x)|}\over {\lambda}} \right)w(x)\,dx
\end{equation}
where $\Phi(t)=t(1+\log^{+}t)$.
\end{thm}

\section{Preliminaries and notation}\label{preliminares}

In this section we gather  some well known  definitions and
properties which we will use along this paper.

\textbf{Maximal operators.} Given a locally integrable function $f$
on $\RN$, the Hardy--Littlewood maximal operator $M$ is defined by
\begin{equation}
Mf(x)=\sup_{\displaystyle{Q\ni x}}{1\over {|Q|}}\int_{Q} f(y)\\,dy,
\end{equation}
where the supremum is taken over all cubes $Q$ containing the point
$x$.

Also we will use  the following operator:
$$M_\delta ^\# f(x)=(M^\# (|f|^\delta
)(x))^{1/\delta}$$
where $M^\#$ is the usual sharp maximal function of C.
Fefferman-Stein:
$$M^\#(f)(x)=\sup_{Q\ni x}\frac{1}{|Q|}\int_Q |f(y)-f_Q|\, dy,  $$
and as usual  $f_Q=\frac{1}{|Q|}\int_Q f(y)\, dy$, and
$$
M_\varepsilon  f(x)=(M (|f|^\varepsilon)(x))^{1/\varepsilon}.
$$

If the supremum is restricted to the dyadic cubes, we will use
respectively the following notation  $M^d$, $M_\delta ^{\#,d}$ and $M_\delta ^{d}$.

\quad

\textbf{Calder\'{o}n--Zygmund operators.} We will use standard well
known definitions, see for instance \cite{J,GrMF}. Let $K(x,y)$ be a
locally integrable function defined of the diagonal $x=y$ in $\RN
\times \RN$, which satisfies the size estimate
\begin{equation}\label{tamano}
|K(x,y)|\leq {c\over{|x-y|^n}},
\end{equation}
and for some $\varepsilon>0$, the regularity condition
\begin{equation}\label{regularidad}
|K(x,y)-K(z,y)|+|K(y,x)-K(y,z)|\leq
c{{|x-z|^{\varepsilon}}\over{|x-y|^{n+\varepsilon}}},
\end{equation}
whenever $2|x-z|<|x-y|$.

A linear operator $T:C_{c}^{\infty}(\RN)\longrightarrow
L_{loc}^{1}(\RN)$ is a Calder\'{o}n--Zygmund operator if it extends
to a bounded operator on $L^2(\RN)$, and there is a kernel $K$
satisfying \eqref{tamano} and \eqref{regularidad} such that
\begin{equation}
Tf(x)=\int_{\RN}K(x,y)f(y)\, dy,
\end{equation}
for any $f\in C_{c}^{\infty}(\RN)$ and $x\notin supp(f)$.

We will use the following result from \cite{APe} several times.

\begin{lemma}\cite{APe}\label{maximalsharp}
Let $T$ be a Calder\'on--Zygmund operator, and $0<\e<1$. Then there
exists a constant $c_{\e}$ such that
\begin{equation}
M_{\varepsilon}^{\#}(Tf)(x)\leq c_{\e}\,Mf(x).
\end{equation}
\end{lemma}

\quad

\textbf{Commutators.} Let $T$ be any operator and let $b$ be any
locally integrable function. The commutator operator $[b, T]$ is
defined by
$$
[b,T]f=b\,T(f)-T(bf).
$$
As we already mentioned when $T$ is any Calder\'{o}n--Zygmund
operator this operator is a bounded operator on $L^{p}(\R^{n})$,
$1<p<\infty$ (\cite{CRoW}) but it is not of weak type $(1,1)$ when $b
\in BMO$ and the following result holds \cite{Pe1}.
\begin{thm}\cite{Pe1} \label{weakcomm}
Let $b$ be a $BMO$ function, $w\, \in A_1$ and $T$ be a singular
integral. Then there exists a positive constant $c=c_{\|b\|_{BMO},
[w]_{A_1}}$ such that for all compact support function $f$ and for
all $\lambda
>0$
$$
w(\{x\in \R^{n}: |[b,T]f(x)|>\lambda\}) \leq c_{\|b\|_{BMO},
[w]_{A_1}}\, \int_{\R^{n}}
\Phi\left(\dfrac{|f(x)|}{\lambda}\right)w(x)\,dx,
$$
being $\Phi(t)= t(1+\log^{+}t)$.
\end{thm}

Later on, G. Pradolini and C. P\'{e}rez in \cite{PePr} improved this
result as follows: given \,$\e>0$
$$
w(\{x\in \R^{n}: |[b, T]f(x)|>\lambda\}) \leq
c\,\Phi\left(\|b\|_{BMO}\right)\,  \int_{\R^{n}}
\Phi\left(\dfrac{|f(x)|}{\lambda}\right) \,M_{\strt{1.7ex} L (\log
L)^{1+\e} }w(x)dx.
$$
where $\Phi(t)= t(1+\log^{+}t)$. The constant $c$ is independent of
the weight $w$, $f$ and $\lambda>0$. The point here is that there is
no condition on the weight $w$.

One of the important steps to prove the last result is to give a version
of the following classical result due to Coifman and C. Fefferman
\cite{CF}: let $T$ be any Calder\'on--Zygmund operator and let
$0<p<\infty$, then there exists a constant $c=c_{p,[w]_{A_\infty}}$
such that for any $w\in A_{\infty}$,
\begin{equation}\label{coifman1}
\int_{\mathbb R^n} |Tf(x)|^p w(x) dx \leq c_{p,[w]_{A_\infty}\,}
\int_{\mathbb R^n} Mf(x)^p w(x)dx.
\end{equation}
The corresponding version for commutators was found in \cite{Pe2}.
\begin{thm}\cite{Pe2}\label{coifman-fefermanConmu}
Let $0 < p < \infty$ and let $w \in A_{\infty}$. Then there exists a
positive constant $c=c_{p,[w]_{A_\infty},\|b\|_{BMO}}$ such that
\begin{equation}\label{conmutadores}
\int_{\R^{n}} \absval{[b,T]f(y)}^{p}\, w(y)\,dy \le
c_{p,[w]_{A_\infty},\|b\|_{BMO}} \, \int_{\R^{n}} M^{2}f(y)^{p}\,
w(y)\,dy.
\end{equation}
\end{thm}

Besides, we will need the following pointwise inequality for
commutators:

\begin{lemma}\cite{Pe1}\label{deltaepsilon}
Let $b\ in BMO$  and let $0 < \delta< \varepsilon$. Then there
exists a positive constant $c= c_{\delta,\varepsilon}$ such that,
\begin{equation}
M_{\delta}^{\#}([b,T]f)(x)\leq c \norm
{b}_{BMO}\left(M_{\varepsilon}(Tf)(x)+M^2f(x)\right),
\end{equation}
for all smooth functions f.
\end{lemma}

\textbf{Orlicz maximal functions.}  We need some few definitions and facts about Orlicz spaces. (For
more information, see Bennett and Sharpley \cite{BS} or Rao and Ren
\cite{RRe}).  A function $B : [0,\infty)\rightarrow [0,\infty)$ is a
doubling Young function if it is continuous, convex and increasing,
if $B(0)=0$ and $B(t)\rightarrow \infty$ as $t\rightarrow\infty$,
and if it satisfies $B(2t)\leq CB(t)$ for all $t>0$.

Recall that we defined the localized Luxembourg norm by equation
$M_{\strt{1.7ex} A}=M_{\strt{1.7ex} A(L)}$, where $M_{\strt{1.7ex} A(L)}$ denotes a maximal type
function defined by the expression
\[
M_{\strt{1.7ex} A(L)}f(x) = \sup_{Q\ni x} \|f\|_{A,Q},
\]
where $A$ is any Young function and $\|f\|_{A,Q}$ denotes the
$A$--average over $Q$ defined by means of the Luxemburg norm
\begin{equation}
\|f\|_{A,Q} = \inf\left\{  \lambda > 0 : \frac{1}{|Q|} \int_Q
A\left(\frac{|f|}{\lambda}\right)\,dx \leq 1 \right\}.
\label{luxembourg}
\end{equation}

An equivalent norm which is often useful in calculations is the following see
Rao and Ren \cite[p.\ 69]{RRe}:
\begin{equation}
\|f\|_{A,Q} \leq  \inf_{\mu>0}\left\{ \mu + \frac{\mu}{|Q|} \int_Q
A\left(\frac{|f|}{\mu}\right)\,dx \right\} \leq 2\|f\|_{A,Q}.
\label{RaoRen}
\end{equation}
Given a Young function $A$, $\bar{A}$ will denote the complementary
Young function associated to $A$; it has the property that for all
$t>0$,
\[ t \leq A^{-1}(t)\bar{A}^{-1}(t) \leq 2t. \]
The property that we will be using is the following generalized
H\"older inequality
\begin{equation}
\frac{1}{|Q|} \int_Q |fg| \leq 2\|f\|_{A,Q}\|g\|_{\bar{A},Q}.
\label{GHI}
\end{equation}

\textbf{B.M.O. functions and John-Nirenberg inequality.} Denote the following pair of complementary Young functions
$$
\Phi(t)= t\,(1+\log^+t) \,\,\, {\rm and}\,\,\, \Psi(t)=e^t-1,
$$
defining the classical Zygmund spaces $L(\log L)$,  and $\exp L$
respectively. The corresponding averages will be denoted by
$$
\|\cdot\|_{\Phi,Q} = \|\cdot\|_{L(\log L),Q}\,\,\, {\rm and}
\,\,\,\|\cdot\|_{\Psi,Q} = \|\cdot\|_{\exp L,Q}.
$$
Using this notation and the generalized H\"older inequality we also get
\begin{equation}\label{Holder-Orlicz2}
\frac1{|Q|}\,\int_{Q}|f(x)\,g(x)|\,dx \le c\,\|f\|_{\exp
L,Q}\,\|g\|_{L(\log L),Q}.
\end{equation}
This inequality allows to write the following formula that will be
used in this article:
\begin{equation}\label{John-Nirenberg-LLogL}
\frac{1}{|Q|}\int_Q \left |b(y)-b_Q \right |f(y)\,dy\leq
c\,\|b\|_{BMO}\|f\|_{L(\log L),Q}.
\end{equation}
for any function $b \in BMO$ and any non negative function $f$. This
inequality follows from \eqref{Holder-Orlicz2} and the
John-Nirenberg inequality \cite{J} for $BMO$ functions: there are
dimensional positive constants $c_1<1$  and $c_2>2$ such that
$$
\frac{1}{  |Q| } \int_{Q} \exp \left( \frac{ c_1 |b(y)- b_{Q}| }{
\|b\|_{BMO}  } \right)\, dy \le c_2
$$
which easily implies that for  appropriate constant $c>0$
$$
\|b-b_Q\|_{\exp L,Q} \le  c\, \|b\|_{BMO}.
$$

\textbf{Weights.} We say that a weight $w$ satisfies the $A_p$
condition for $1<p<\infty$ if there exists a constant $c$ such that
for any cube $Q$,
\begin{equation*}
\left( {\displaystyle{1\over{|Q|}}}\int_{Q}w \right)
\left({\displaystyle{1\over{|Q|}}}\int_{Q}w
^{1-p'}\right)^{p-1}=c<\infty.
\end{equation*}
We will denote as $[w]_{A_p}$ the smallest of these $c$.

Also we recall that $w$ is an $A_1$ weight if there is a finite
constant $c$ such that $Mw\le c\,w$ a.e., and where $[w]_{A_1}$
denotes the smallest of these $c$. Also as usual, we denote
$A_{\infty}=\bigcup_{\displaystyle{p\geq 1}}A_p$.

We will use several well-know properties about the $A_p$ weights.
First, it follows from H\"{o}lder's inequality that if $w_1, w_2\in
A_1$, then $w=w_1w_2^{1-p}\in A_p$, and
\begin{equation}\label{pro1peso}
[w]_{A_p}\leq [w_1]_{A_1}[w_2]_{A_1}^{p-1}.
\end{equation}

Second, if $0< \delta <1$, then $(Mf)^{\delta}\in A_{1}$ (see
\cite{CRo}), and $f\in L^1_{loc}(\RN)$
\begin{equation}\label{CRo}
[(Mf)^{\delta}]_{A_1}\leq {\displaystyle{c_n\over{1-\delta}}}
\end{equation}

Third, is the Fefferman-Stein inequality \cite{FS} saying that for
any weight $w$,
\begin{equation}\label{fs1}
\|Mf\|_{\strt{1.7ex}L^p(w)}\le c_n\, p'\,
\|f\|_{\strt{1.7ex}L^p(Mw)}\quad (1<p<\infty),
\end{equation}
where, as usual, $p'$ denotes the dual exponent of $p$,
$p'=\frac{p}{p-1}$.

\textbf{Sharp reverse H\"{o}lder's inequality.}  An important step toward the proof of Theorem \ref{fuerte} is the following lemma which gives a precise version of the following classical reverse H\"{o}lder's inequality: if $w\in A_1$, there are
constants $r>1$ and $c\geq 1$ such that
\begin{equation}
M_{r}w(x)\leq cw(x).
\end{equation}

\begin{lemma}\cite{LOPe1}\label{roptimo}
Let $w\in A_1$, and let $r_{w}=1+{1\over{2^{n+1}[w]_{A_1}}}$. Then,
\begin{equation}
M_{r_{w}}w(x)\leq 2[w]_{A_1}w(x).
\end{equation}
\end{lemma}

\textbf{Building $A_1$ weights from duality.}  The following lemma gives a way to produce $A_{1}$ weights with special control on the constant. It is based on the  so called Rubio de Francia iteration scheme or algorithm.

\begin{lemma}\cite{LOPe3}\label{rubiodefrancia}
Let $1<s<\infty$, and let $v$ be a weight. Then, there exists a
nonnegative sublinear operator $R$ bounded in $L^s(v)$ satisfying the following
properties:
\begin{enumerate}
  \item[(i)] $h\leq R(h)$;
  \item[(ii)] $\norm{Rh}_{L^{s}(v)}\leq 2\norm{h}_{L^{s}(v)}$;
  \item[(iii)] ${Rh\,v^{1/s}\in A_1}$ with
  \begin{equation}
 [Rh\,v^{1/s}]_{A_1}\leq cs'.
  \end{equation}
\end{enumerate}
\end{lemma}

\textbf{The sharp maximal function with optimal constants.}  We will need the following result linking the weighted $L^{p}$ norm of a function and its sharp maximal function. This result is central in our approach. It can be found in \cite{Pe4}.

\begin{lemma}\cite{Pe4}\label{bagby}
Let $0<p<\infty$, $0<\delta<1$, and let $w\in A_q$, with $1\leq
q<\infty$. Then
there exists a constant $c=c(n,q,\delta)$ such that
\begin{equation}
\norm{f}_{L^{p}(w)}\leq
cp[w]_{A_q}\norm{M_{\delta}^{\#,d}(f)}_{L^{p}(w)},
\end{equation}
for any function $f$ such that $|\{x:|f(x)|>t\}|<\infty$.
\end{lemma}

\textbf{Buckley's Theorem.}  The following well known result gives the sharp bound for the Hardy-Littlewood maximal function on $L^p(w)$.

\begin{thm}\cite{Bu} \label{buck}
Let $M$ be the Hardy--Littlewood maximal function, $p>1$ and $w\in A_p$. Then
there exists a constant $c=c_{n}p'$ such that
\begin{equation}\label{buckin}
\norm{M}_{L^p(w)}\leq c_{n}p'\,[w]_{A_p}^{1\over {p-1}}.
\end{equation}
\end{thm}

S.M. Buckley did not get the constant $c_{n}p'$ that appears in \eqref{buckin} but, recently a simple and elegant proof of this theorem was given by Lerner \cite{L1} where it was obtained.

Observe that we could simply use Theorem \ref{buck}, inequality \eqref{pro1peso} and then
\eqref{CRo} to show that $\norm{M}_{L^p((M_r w)^{1-p})}$ is
bounded by a multiple of $[(M_r w)^{1-p}]_{A_p}^{\frac{1}{p-1}} \leq
[M_r w]_{A_1}^{\frac{p-1}{p-1}} \leq c\,r'$. However this estimate
is worse than the one obtained in the following lemma.

\begin{lemma}\label{chicos} Let $1<p,r<\infty$, then
\begin{equation}\label{chicos2}
\norm{Mf}_{L^p((M_r w)^{1-p})}\leq
c\,p'(r')^{\frac1p}\norm{f}_{L^p( w^{1-p})}.
\end{equation}
\end{lemma}

\noindent {\sc Proof of Lemma \ref{chicos}:}
This is a consequence of Lemma 3.4 from \cite{LOPe1}. Indeed, this
lemma yields
$$
\norm{Mf}_{L^p((M_r w)^{1-p})}\leq
c\,p'(r')^{1-\frac{1}{rp'}}\norm{f}_{L^p( w^{1-p})},
$$
but $1-\frac{1}{rp'} =
\frac{1}{p'}+\frac{1}{p}-\frac{1}{rp'}=\frac{1}{p}+\frac{1}{r'p'}$
and hence $(r')^{1-\frac{1}{rp'}}=(r')^{\frac{1}{p}+\frac{1}{r'p'}}
\leq c\,(r')^{\frac{1}{p}}$ since $t^{1/t}\le e, \quad t\geq 1$.

\hfill $\Box$

\textbf{A two weight inequality for Singular Integrals with optimal bounds.}  The last  important lemma we will use can be seen as a dual version of Lemma \ref{chicos} for Singular Integrals.
It can be found essentially in \cite{LOPe3} (see also \cite{Pe4}).

\begin{lemma}\cite{LOPe3}  \label{sharp1} Let $T$ be a Calder\'on--Zygmund singular integral operator, and let $1<p,r<\infty$. Then there is a constant $c=C_{n,T}$ such that
\begin{equation*}
\|Tf\|_{\strt{1.7ex}L^p(w)}\le cp'\,(r')
^{1/p'}\|f\|_{\strt{1.7ex}L^p(M_rw)} \quad w\geq 0.
\end{equation*}

\end{lemma}

As a corollary of this lemma we have the following estimate which was crucial to derive the main result from \cite{LOPe3} and will be used in the proof of Theorem \ref{debil}.

\begin{coroll}\cite{LOPe3} \label{LOPe3-weakwithparameters}  Let $T$ be a Calder\'on--Zygmund singular integral operator, and let $1<p,r<\infty$. Then there is a constant $c=C_{n,T}$ such that
\begin{equation}\label{casidebil}
 \|Tf\|_{L^{1,\infty}(w)}  \leq c(p')^{p}(r')^{p-1}\int_{\R^n}f(x)\, M_{r}w(x)dx.
\end{equation}
\end{coroll}

\textbf{Kolmogorov inequality.} Finally, we will employ several times
the well known Kolmogorov inequality.  Let $0<p<q<\infty$,
then there is a constant $C=C_{p,q}$ such that for any measurable
function $f$
\begin{equation}\label{kolmogorov}
\|f\|_{L^p(Q, \frac{dx}{|Q|})}\leq C\, \|f\|_{L^{q,\infty}(Q,
\frac{dx}{|Q|})}.
\end{equation}
See for instance  \cite{GrCF} p. 91, ex. 2.1.5.

\section{The strong case}

Before proving Theorem \ref{fuerte}, we need the following two
results. The first one is inspired by Lemma \ref{bagby}, and the
second one is about the $L^p$ boundeness of $M^2$ from
$L^p(w^{1-p})$ to $L^p((M_{r}w)^{1-p})$:

\begin{lemma}\label{mbagby}
Let $0<p<\infty$, $1\leq q<\infty$, $0< \varepsilon \le 1$ and $w\in A_q$. Suppose that $f$
is such that for each $t>0$, $|\{x:f(x)>t\}|<\infty$.  Then there is a constant $c=c_{n,q,\varepsilon}$ such that
\begin{equation}
\norm{M^d_{\varepsilon}f}_{L^{p}(w)}\leq
c\,p[w]_{A_q}\norm{M^{\#,d}_{\varepsilon}f}_{L^{p}(w)}.
\end{equation}
\end{lemma}

\noindent {\sc Proof:}

(In this proof, and for simplicity of  notation we denote $M=M^d$
and similarly for the other operators).

In order to prove Lemma \ref{mbagby}, we apply Lemma \ref{bagby} to $M_{\varepsilon}f$ with $\delta=\varepsilon_{0}$, such that $0<\varepsilon_{0}<\varepsilon<1$, and get
\begin{equation}\label{ineqMsharp}
\norm{M_{\varepsilon}f}_{L^{p}(w)}\leq
cp[w]_{A_q}\norm{M^{\#}_{\varepsilon_{0}}(M_{\varepsilon}f)}_{L^{p}(w)}.
\end{equation}

Then, we will finish if we see that if $0 < \varepsilon_{0} < \varepsilon < 1$
\begin{equation}
M^{\#}_{\varepsilon_{0}}(M_{\varepsilon}f)(x)\leq c\, M^{\#}_{\varepsilon}f(x),
\end{equation}
where recall that if $f\geq 0$
\begin{equation}
M^{\#}_{\varepsilon_{0}}f(x)=M^{\#}(f^{\varepsilon_{0}})^{1/\varepsilon_{0}}(x)=\sup_{\displaystyle{Q\ni
x}}\left({1\over {|Q|}}\int_{Q}
|f^{\varepsilon_{0}}-(f^{\varepsilon_{0}})_Q|\,dy,\right)^{1/\varepsilon_{0}}.
\end{equation}

Now fix $x$  and a dyadic cube $Q$ with $x\in Q$. Hence
\begin{eqnarray}\label{sharp}
{1\over {|Q|}}\int_{Q}\left|
(M(f^{\varepsilon}))^{\varepsilon_{0}/\varepsilon}(y)-((M(f^{\varepsilon}))^{\varepsilon_{0}/\varepsilon})_Q\right|\,dy\\
\leq 2 {1\over {|Q|}}\int_{Q}
(M(f^{\varepsilon}))^{\varepsilon_{0}/\varepsilon}(y)\,dy-
\inf_{\displaystyle{Q}}(M(f^{\varepsilon}))^{\varepsilon_{0}/\varepsilon}
\end{eqnarray}
adding and susbtracting
$\inf_{\displaystyle{Q}}(M(f^{\varepsilon}))^{\varepsilon_{0}/\varepsilon}$.

Now we have,
\begin{equation}
f^{\varepsilon}(x)= g(x)+h(x),
\end{equation}
where
$g(x)=(f^{\varepsilon}(x)-f^{\varepsilon}_{Q})\chi_{\strt{1.7ex}Q}(x)$
and $h(x)=
f^{\varepsilon}_{Q}\chi_{\strt{1.7ex}Q}(x)+f^{\varepsilon}(x)\chi_{\strt{1.7ex}(Q)^c}(x)$.

Then, since $\varepsilon_{0}/\varepsilon<1$,
\begin{equation}
{1\over {|Q|}}\int_{Q} (M(f^{\varepsilon}))^{\varepsilon_{0}/\varepsilon}(y)\,dy
\leq {1\over{|Q|}}\int_{Q} (Mg)^{\varepsilon_{0}/\varepsilon}(y)\, dy
+{1\over{|Q|}}\int_{Q} (Mh)^{\varepsilon_{0}/\varepsilon}(y)\, dy
\end{equation}
To finish, we study each separately. For the first one we use
Kolmogorov's inequality \eqref{kolmogorov} with $p=\varepsilon_{0}/\varepsilon<1=q$, and
the fact that $M$ is of weak-type $(1,1)$,

\begin{eqnarray}
{1\over{|Q|}}\int_{Q} (Mg)^{\varepsilon_{0}/\varepsilon}(y)\, dy &\leq&
C_{\varepsilon, \varepsilon_{0}}\,\left( \frac{1}{|Q|} \int_{Q} |g(y)|\,
dy\right)^{\varepsilon_{0}/\varepsilon}\\
&=&C_{\varepsilon, \varepsilon_{0}} \left(\frac{1}{|Q|}\int_{Q}
|f^{\varepsilon}(y)-f^{\varepsilon}_{Q}|\,
dy\right)^{\varepsilon_{0}/\varepsilon}\\\label{sharpg}
&\leq& C_{\varepsilon, \varepsilon_{0}}\,(M^{\#}_{\varepsilon}f(x))^{\varepsilon_{0}}.
\end{eqnarray}
This part is the bad term because the other term is less singular.
Indeed,  we claim the following
$$
{1\over{|Q|}}\int_{Q} (Mh)^{\varepsilon_{0}/\varepsilon}(y)\, dy \leq
(\displaystyle{\inf_{Q} M(f^\varepsilon)})^{\varepsilon_{0}/\varepsilon}.
$$
Combining this inequality together with \eqref{sharpg} and
\eqref{sharp} we derive \eqref{ineqMsharp} concluding the proof of
the lemma.

To prove the claim we recall $Mh(x)={\displaystyle\sup_{R\ni
x}}{1\over {|R|}}\int_{R} h(y)\,dy$, where the supremum is taken
over any dyadic cube $R$ containing $x$, and we distinguish two types
of cubes:

\begin{enumerate}
  \item let $R\subset Q$.
    In this case,
    \begin{equation}
    {1\over {|R|}}\int_{R} h(y)\,dy={1\over {|R|}}\int_{R}
    f^{\varepsilon}_Q\,dy=f^{\varepsilon}_Q\leq \displaystyle \inf_{Q}M(f^{\varepsilon}).
    \end{equation}
  \item $R\supset Q$
    In this case
    \begin{eqnarray}
    {1\over {|R|}}\int_{R} h(y)\,dy&=&{1\over {|R|}}|R\cap Q| f^{\varepsilon}_Q+{1\over {|R|}}\int_{R\cap Q^c}f^{\varepsilon}(y)\,dy \nonumber\\
    &=&{{|Q|}\over {|R|}}f^{\varepsilon}_Q+{1\over {|R|}}\int_{R\cap
    Q^c}f^{\varepsilon}(y)\,dy\nonumber\\
    &=&{1\over {|R|}}\int_{R\cap Q} f^{\varepsilon}(x)\,dx+{1\over {|R|}}\int_{R\cap Q^c}
    f^{\varepsilon}(x)\,dx\nonumber\\
    &=&{1\over {|R|}}\int_{R} f^{\varepsilon}(x)\,dx\nonumber \leq \displaystyle{\inf_{Q}} M(f^{\varepsilon}).
    \end{eqnarray}
\end{enumerate}
So,
\begin{equation}
{1\over{|Q|}}\int_{Q} (Mh)^{\varepsilon_{0}/\varepsilon}(y)\, dy \leq
{1\over{|Q|}}\int_{Q} (\inf_Q M(f^{\varepsilon}))^{\varepsilon_{0}/\varepsilon}\, dy =
(\displaystyle{\inf_{Q}M(f^{\varepsilon})})^{\varepsilon_{0}/\varepsilon}.
\end{equation}

\hfill $\Box$

\begin{prop}\label{m2dospesos}
Let $M^2$ be the composition $M\circ M$, and let $1<p,r<\infty$.
Then, there is a constant $c$ independent of $r,p$ such that
$$
\norm{M^2f}_{L^p((M_{r}w)^{1-p})}\le c\,(p')^2(r')^{1+1/p}
\norm{f}_{L^p(w^{1-p})}.
$$

\end{prop}

\noindent {\sc Proof:}

To prove the inequality we use Buckley's Theorem \ref{buck}
\begin{eqnarray*}
\norm{M^{2}f}_{L^p((M_{r}w)^{1-p})}&=&\norm{M(Mf)}_{L^p((M_{r}w)^{1-p})}\\
&\leq&c_{n}p'\,[(M_{r}w)^{1-p}]_{A_p}^{1/(p-1)}\norm{Mf}_{L^p((M_{r}w)^{1-p})}\\
&\leq& c_{n}p'\,[(M_{r}w)]_{A_1}^{(p-1)/(p-1)}\norm{Mf}_{L^p((M_{r}w)^{1-p})}\\
&\leq&c_{n}\,r'(p')^2(r')^{1/p}\norm{f}_{L^p(w^{1-p})}
\end{eqnarray*}
where we have used property \eqref{pro1peso} and Lemma \ref{chicos}.

\hfill $\Box$

\noindent {\sc Proof of Theorem \ref{fuerte}:} We will prove
\eqref{acotadospesos}, namely
\begin{equation}\label{acotadospesos2}
\norm{[b,T]f}_{L^{p}(w)}\le c {(pp')}^2\, \norm{b}_{BMO}\,
(r')^{1+\frac{1}{p'}}\,
   \norm{f}_{L^{p}(M_{r}w)}.
\end{equation}
A direct application of Lemma \ref{roptimo} with
$r=r_{w}=1+{1\over{2^{n+1}[w]_{A_1}}}$  would finish the proof of
the Theorem \ref{fuerte}.

By duality \eqref{acotadospesos2} is equivalent to proving
\begin{equation}\label{acotadospesosdual}
\norm{\displaystyle{{[b,T]^*f}\over{M_{r}w}}}_{L^{p'}(M_{r}w)}\le
c\, (pp')^2\, \norm{b}_{BMO}\,
(r')^{1+\frac{1}{p'}}\,\norm{\displaystyle{{f}\over{w}}}_{L^{p'}(w)},
\end{equation}
where $[b,T]^{*}$ is the adjoint operator of $[b,T]$ with respect to the $L^2-$pairing. Now,

\begin{equation}
\norm{\displaystyle{{[b,T]^{*}f}\over{M_{r}w}}}_{L^{p'}(M_{r}w)}=\sup_{\displaystyle{\norm{h}_{L^{p}(M_{r}w)}=1}}
\left|\int_{\RN}[b,T]^{*}f(x)h(x)\, dx \right|.
\end{equation}

First by part (i) of Lemma \ref{rubiodefrancia} for $s=p'$ \, and \, $v=M_rw$\, there
exists an operator~$R$ such that
\begin{eqnarray*}
 I&=&\left|\int_{\RN}[b,T]^{*}f(x)h(x)\, dx\right|\\
 &\leq&\int_{\RN}|[b,T]^{*}f(x)|\,|h(x)|\, dx \\
 &\leq& \int_{\RN}|[b,T]^{*}f(x)|\,Rh(x)\, dx.
\end{eqnarray*}

Besides, combining part (iii) of Lemma \ref{rubiodefrancia} for $s=p$ and \, $v=M_rw$  with \eqref{pro1peso} and \eqref{CRo} we obtain
 \begin{eqnarray}
[Rh]_{A_3}&=&[Rh (M_rw)^{1/p}(M_rw)^{-1/p}]_{A_3}= [Rh (M_rw)^{1/p}((M_rw)^{1/2p})^{-2}]_{A_3}\\
&\leq& [Rh (M_rw)^{1/p}]_{A_1}[(M_rw)^{1/2p}]_{A_1}^2\\
\label{a3}&\leq& c_{n}\,p'.
\end{eqnarray}

Applying now Lemma \ref{bagby} to $[b,T]^{*}$ with weight $w=Rh$, $q=3$ and $p=1$ we get
\begin{eqnarray*}
I\leq
c_{\delta}\,[Rh]_{A_3}\int_{\RN}M_{\delta}^{\#}([b,T]^{*}f)(x)\,Rh(x)\,dx.
\end{eqnarray*}

Then we apply Lemma \ref{deltaepsilon} with
$0<\d<\e$ to $[b,T]^{*}=-[b,T^*]$ since is also a commutator with a
Calder\'{o}n-Zygmunnd operator.
By H\"{o}lder's Inequality and using the property (ii) of the Lemma
\ref{rubiodefrancia} and \eqref{a3} we can continue with
\begin{eqnarray*}
I&\leq&
c_{n,\delta,\varepsilon}\norm{b}_{BMO}p'\int_{\RN}\left(M_{\varepsilon}(T^*f)(x)+M^2f(x)\right)\,Rh(x)\,
dx\\
&=& c_{n,\delta,\varepsilon}\norm{b}_{BMO}p'(I_{1}+I_{2}).
\end{eqnarray*}
Now
\begin{eqnarray*}
I_2&=&\int_{\RN}M^2f(x)\,Rh(x)\,dx \\
&\leq& \left(\int_{\RN}|M^{2}f(x)|^{p'}(M_{r}w(x))^{1-p'}\,dx\right)^{1/p'}\left(\int_{\RN}M_rw(x)(Rh(x))^p\,dx\right)^{1/p}\\
&=& 2\,\norm{\displaystyle{{M^2f}\over{M_{r}w}}}_{L^{p'}(M_{r}w)}.
\end{eqnarray*}

For the first term $I_1$, we apply Lemma \ref{mbagby} to
$M_\varepsilon(T^*f)$ with weight $w=Rh$, $q=3$ and $p=1$, and then
Lemma \ref{maximalsharp} (choosing now $0<\varepsilon<1$):
\begin{eqnarray*}
I_1&=&\int_{\RN}|M_{\varepsilon}(T^*f)(x)|\,|Rh(x)|dx\leq
c_{n,\varepsilon}\,[Rh]_{A_3}\int_{\RN}|M_{\varepsilon}^{\#}(T^*f)(x)|\,|Rh(x)|dx\\
&\leq&c_{n,\varepsilon}\,[Rh]_{A_3}\,\int_{\RN}|Mf(x)|\,|Rh(x)|dx\leq c_{n,\varepsilon}\, p'\,
\left(\int_{\RN}Mf(x)^{p'}(M_{r}w(x))^{1-p'}\,dx\right)^{1/p'}\\
&=&c_{n,\varepsilon}\, p'\norm{\displaystyle{{Mf}\over{M_{r}w}}}_{L^{p'}(M_{r}w)},
\end{eqnarray*}
using property (ii) of Lemma \ref{rubiodefrancia} and H\"{o}lder's
inequality, as we did when estimating $I_2$.

Therefore, combining estimates we have
\begin{equation}
\norm{\displaystyle{{[b,T]^{*}}\over{M_{r}w}}}_{L^{p'}(M_{r}w)}\leq
c_{n,\delta,\varepsilon}\,\norm{b}_{BMO}\,(p')^2\norm{M^2f}_{L^{p'}((M_{r}w)^{1-p'})}.
\end{equation}
Finally, to finish the proof of the theorem we apply Proposition \ref{m2dospesos} to get the estimate \eqref{acotadospesosdual}  which yields the claim
\eqref{acotadospesos2}.

\hfill $\Box$

\section{The weak case}

In this section we will prove Theorem \ref{debil}. For $f\in C^\infty_0(\mathbb R^n)$ we consider the classical
Calder\'on-Zygmund decomposition of $f$ at level $\lambda$.
Therefore we get a family $\{Q_j\}=Q_{j}(x_{Q_{j}},r_{j})$ of
non-overlapping dyadic cubes satisfying

\begin{equation}\label{CZdec}
\lambda<\frac{1}{|Q_j|}\int_{Q_j}|f(x)|\,dx\leq 2^n\lambda.
\end{equation}

This implies that for $\Omega=\displaystyle{\cup_j} Q_j$ we have
$|f(x)|\leq \lambda$ a.e. on $\mathbb R^n\setminus \Omega$.

Now we split $f$ in the standard ``good'' and ``bad'' parts $f=g+h$.
Indeed, if as usual $f_{Q_j}$ denotes the average of $f$ on $Q_j$,
we take
$$
g(x)=\left \{ \begin{array}{ll} f(x),&x\in \mathbb R^n\setminus \Omega\\
f_{Q_j},&x\in Q_j,\end{array}\right.
$$
which also verifies $|g(x)|\leq 2^n\lambda$ a.e. For the bad part we
consider $h=\displaystyle{\sum_j h_j}$ where
$h_j(x)=(f(x)-f_{Q_j})\chi_{Q_j}(x)$.

For each $j$ we denote
$$\widetilde{Q_j}=3Q_j, \quad \widetilde{\Omega}=\displaystyle{\cup_j \widetilde{Q_j}}\quad  \mbox{and} \quad
w_{j}(x)=w(x)\chi_{_{ \R^{n}\setminus 3Q_{j}}}$$

Then we can write
\begin{eqnarray*}
w(\{x\in \RN:|[b,T]f(x)|>\lambda\})&\leq&w(\{x\in \RN \setminus \widetilde{\Omega} :|[b,T]g(x)|>\lambda/2\})\\
&+&w(\widetilde{\Omega})\\
&+&w(\{x\in \RN\setminus
\widetilde{\Omega}:|[b,T]h(x)|>\lambda/2\})\\
&=&I+II+III.
\end{eqnarray*}
and we study each term separately.

For the first one we use Chebyschev's inequality, calling
$\widetilde{w}(x)=w(x)\chi_{\RN \setminus \widetilde{\Omega}}(x)$,
and by the $L^p$ estimate \eqref{acotadospesos} from Theorem
\ref{fuerte} which holds for any weight we have
\begin{eqnarray*}
I &\leq& \dfrac{c}{\lambda^{p}}\int_{\R^{n}}|[b,T]g(x)|^{p}
\widetilde{w}(x)\, dx\\
&\leq& \dfrac{c}{\lambda^{p}}\|b\|_{BMO}^p(pp')^{2p}
(r')^{(1+{1\over {p'}})p}\int_{\R^{n}}|g(x)|^{p}
\,M_{r}\widetilde{w}(x)\,dx\\
&\leq& \dfrac{c}{\lambda}\|b\|_{BMO}^p(pp')^{2p}(r')^{
2p-1}\int_{\R^{n}}|g(x)|
\,M_{r}\widetilde{w}(x)\,dx\\
&=& \dfrac{c}{\lambda}\|b\|_{BMO}^p(pp')^{2p}(r')^{2p-1}\left(
\int_{ \R^n \setminus \Omega }|f(x)| \,M_{r}\widetilde{w}(x)\,dx+
\int_{\Omega}|g(x)| \,M_{r}\widetilde{w}(x)\,dx \right).\\
\end{eqnarray*}
It is clear that we only need to estimate the second term in the
last expression and to do so we will be using the following fact:
For any nonnegative function $u$ with $M u(x)< \infty $ a.e., any
cube $Q$, and any $R>1$ we have
\begin{equation} \label{aboutM}
M( \chi_{ \R^{n} \setminus RQ } u )(y) \approx M( \chi_{ \R^{n}
\setminus RQ } u )(z)  \qquad y,z \in Q
\end{equation}
with dimensional constants. This can be found in \cite{GCRdF} p. 159.

Hence we can continue estimating the second term with
\begin{eqnarray*}
\int_{\Omega}|g(x)| \,M_{r}\widetilde{w}(x)\,dx &\leq&
\sum_{j}\int_{Q_{j}}|f_{Q_{j}}|
\, M_r(w_{j})(x)\,dx\\
&=&\sum_{j}\left(\int_{Q_{j}}|f(x)|dx\right)
\dfrac{1}{|Q_j|}\int_{Q_{j}} \,M_{r}w_{j}(x)\,dx\\
&\leq& c\, \sum_{j}\left(\int_{Q_{j}}|f(x)|\,dx\right)
\inf_{x \in Q_{j}}
M_{r}w_{j}\\
&\leq& c\,\int_{\R^{n}}|f(x)| \,M_{r}w(x)\, dx.
\end{eqnarray*}

Up to now the parameter $1<r<\infty$ was arbitrary but if we choose
now $r=r_{w}$ where $r_w$ is the sharp reverse H\"older exponent
$r_{w}=1+{1\over{2^{n+1}[w]_{A_1}}}$ from Lemma \ref{roptimo} we can
continue with
\begin{eqnarray*}
I&\leq&{c\over \lambda} \|b\|_{BMO}^p(pp')^{2p}[w]_{A_1}^{2p-1}
\int_{\R^n }|f(x)| \,M_{r}w(x)\,dx \\
&\leq&{c\over \lambda}\|b\|_{BMO}^p(pp')^{2p}[w]_{A_1}^{2p}
\int_{\R^n }|f(x)|w(x)\,dx.
\end{eqnarray*}

The second term, $II$, is estimated in a very standard way
\begin{eqnarray}\label{segundo}
II=w(\widetilde{\Omega})&\leq &c\,\sum_j \frac{w(\widetilde{Q_j})}
{|\widetilde{Q_j}|}|Q_j|\nonumber \\
&\leq &\frac{c}{\lambda} \sum_j
\frac{w(\widetilde{Q_j})}{|\widetilde{Q_j}|}
 \int_{Q_j} |f(x)|\,dx\nonumber \\
&\leq &\frac{c}{\lambda} \sum_j \int_{Q_j} |f(x)|Mw(x)\,dx\nonumber  \\
&\leq &\frac{c}{\lambda}\int_{\mathbb R^n} |f(x)|Mw(x)\,dx\\
&\leq&\frac{c}{\lambda}[w]_{A_1}\int_{\mathbb R^n} |f(x)|w(x)\,dx,
\end{eqnarray}

because $w\in A_1$.

For part $III$ first note that
$$ [b,T]h(x)=\sum_{j}[b,T]h_{j}(x)=
\sum_{j}(b(x)-b_{Q_{j}})Th_{j}(x) - \sum_{j}
T((b-b_{Q_{j}})h_{j})(x)
$$
where as before $b_{Q}=\frac{1}{|Q|}\int_{Q}b$. Then
\begin{eqnarray*}
III &\leq& w(\{x\in \R^{n} \setminus \widetilde{\Omega}:
|\sum_{j}(b(x)-b_{Q_{j}})Th_{j}(x)|>\dfrac{\lambda}{4}\})\\\nonumber
&&+w(\{x\in \R^{n}\setminus \widetilde{\Omega}:
|\sum_{j}T((b-b_{Q_{j}})h_{j})(x)|>\dfrac{\lambda}{4}\})\\
\nonumber &=&A + B. \nonumber
\end{eqnarray*}

Using the standard estimates of the kernel $K$ we get

\begin{eqnarray*}
A &\leq& \dfrac{c}{\lambda}\int_{\R^{n}\setminus
\widetilde{\Omega}}\sum_{j}|b(x)-b_{Q_{j}}||Th_{j}(x)|\widetilde{w}(x) \,dx \\
&\leq&\dfrac{c}{\lambda}\sum_{j}\int_{\R^{n}\setminus 3Q_j
}|b(x)-b_{Q_{j}}|w_{j}(x)\int_{Q_{j}}|h_{j}(y)||K(x,y)-K(x,x_{Q_{j}})|
\,dy\, dx \\
&\leq&\dfrac{c}{\lambda}\sum_{j}\int_{Q_{j}}|h_{j}(y)|\int_{\R^{n}\setminus
3Q_j}|K(x,y)-K(x,x_{Q_{j}})||b(x)-b_{Q_{j}}|w_{j}(x)
\, dx\, dy \\
&\leq&\dfrac{c}{\lambda}\sum_{j}\int_{Q_{j}}|h_{j}(y)|
\sum_{k=1}^{\infty}\int_{2^{k}r_{j}\leq|x-x_{Q_{j}}|<2^{k}r_{j}}\dfrac{|y-x_{Q_{j}}|^{\varepsilon}}{|x-x_{Q_{j}}|^{n+\varepsilon}}|b(x)-b_{Q_{j}}|w_{j}(x)
\, dx\, dy \\
&\leq&\dfrac{c}{\lambda}\sum_{j}\int_{Q_{j}}|h_{j}(y)|
\left(\sum_{k=1}^{\infty}\dfrac{2^{-k\varepsilon}}{(2^{k}r_{j})^{n}}
\int_{|x-x_{Q_{j}}|<2^{k+1}r_{j}}|b(x)-b_{Q_{j}}|w_{j}(x) \,
dx\right)\,dy
\end{eqnarray*}
To control the sum on $k$ we use standard estimates together with
the generalized H\"older inequality and John-Nirenberg's Theorem.
Indeed if $y \in Q_j$ we have
\[
\sum_{k=1}^{\infty}\dfrac{2^{-k \varepsilon}}{(2^{k+1}r_{j})^{n}}
\int\limits_{|x-x_{ Q_{j} }|<2^{k+1}r_{j}}|b(x)-b_{Q_{j}}|w_{j}(x) \, dx
\qquad \qquad \qquad \qquad \qquad \qquad \qquad
\]
\begin{eqnarray*}
&\le& c\,
\sum_{k=1}^{\infty}\dfrac{2^{-k \varepsilon}}{(2^{k+1}r_{j})^{n}}\int_{2^{k+1}Q_{j}}|b(x)-b_{2^{k+1}Q_{j}}|w_{j}(x)
\, dx \\
&&+\sum_{k=1}^{\infty}\dfrac{2^{-k \varepsilon}}{(2^{k+1}r_{j})^{n}}\int_{2^{k+1}Q_{j}}|b_{2^{k+1}Q_{j}}-b_{Q_{j}}|w_{j}(x)
\, dx\\
&\le& c\, \sum_{k=1}^{\infty}2^{-k \varepsilon}
\|b-b_{2^{k+1}Q_{j}}\|_{\strt{1.7ex}\exp
L,2^{k+1}Q_{j}}\|w_{j}\|_{\strt{1.7ex}L\log L,2^{k+1}Q_{j}}\\
&&+\sum_{k=1}^{\infty}2^{-k \varepsilon}(k+1)\|b\|_{BMO}\,M w_{j}(y)\\
&\le& c\, \sum_{k=1}^{\infty}2^{-k \varepsilon} \|b\|_{BMO}M_{LlogL}w_{j}(y)
+\sum_{k=1}^{\infty}2^{-k \varepsilon}(k+1)\|b\|_{BMO}\,[w]_{A_1}w_{j}(y)\\
&\le& c\, \sum_{k=1}^{\infty} 2^{-k \varepsilon}(k+1)
\|b\|_{BMO}[w]_{A_1}M_{LlogL}w_{j}(y).
\end{eqnarray*}

Then we can continue the estimate of $A$ as follows
\begin{eqnarray*}
A &\leq&\dfrac{c}{\lambda}\sum_{j}\int_{Q_{j}}|h_{j}(y)|\|b\|_{BMO}[w]_{A_1}\,M_{LlogL}w_{j}(y)dy\\
&\leq&\dfrac{c}{\lambda}|\|b\|_{BMO}[w]_{A_1}\sum_{j}\left(\int_{Q_j}|f(y)|\,M_{LlogL}w_{j}(y)dy +\int_{Q_j}|f_{Q_j}|\,M_{LlogL}w_{j}(y)dy\right)\\
&\leq&\dfrac{c}{\lambda}|\|b\|_{BMO}[w]_{A_1}\left(\int_{\RN}|f(y)|\,M_{LlogL}w(y)dy+\sum_{j}\int_{Q_j}|f(x)|\,dx{1\over {|Q_{j}|}}\int_{Q_{j}}M_{LlogL}w_{j}(y)dy\right)\\
&\leq&\dfrac{c}{\lambda}|\|b\|_{BMO}[w]_{A_1}\left(\int_{\RN}|f(y)|\,M_{LlogL}w(y)dy+\sum_{j}\int_{Q_j}|f(x)|\inf_{\displaystyle{Q_{j}}}M_{LlogL}w_{j}\right)\\
&\leq&\dfrac{c}{\lambda}|\|b\|_{BMO}[w]_{A_1}\int_{\RN}|f(y)|\,M_{LlogL}w(y)dy
\end{eqnarray*}

To estimate $B$ we will use inequality   \eqref{casidebil} from Corollary \ref{LOPe3-weakwithparameters} as follows
\begin{eqnarray*}
B &\leq&\widetilde{w}
(\{x\in\R^{n}:|T(\sum_{j}(b-b_{Q_{j}})h_{j})(x)|>\dfrac{\lambda}{4}\})\\
&\leq&c\,{{(p')^p(r')^{p-1}}\over{\lambda}}\int_{\R^{n}}\left|\sum_{j}(b(x)-b_{Q_{j}})h_{j})(x)\right|\,M_{r}(\widetilde{w})(x)dx\\
&\leq&c\,{{(p')^p(r')^{p-1}}\over{\lambda}}\sum_{j}\int_{Q_{j}}|b(x)-b_{Q_{j}}||f(x)-f_{Q_{j}}|
\,M_{r}w_{j}(x)dx\\
&\leq&c\,{{(p')^p(r')^{p-1}}\over{\lambda}}\sum_{j}
\left(\int_{Q_{j}}|b(x)-b_{Q_{j}}||f(x)|\,M_{r}w_{j}(x)dx+
\int_{Q_{j}}|b(x)-b_{Q_{j}}||f_{Q_{j}}|\, M_{r}w_{j}(x)dx\right)\\
&=&c\,{(p')^p(r')^{p-1}}(B_1+B_2) \leq c\,(p')^{p}[w]_{A_{1}}^{p-1}(B_1+B_2)
\end{eqnarray*}

The estimate for $B_2$ we use    \eqref{aboutM}
\begin{eqnarray*}
B_2 &=&\dfrac{c}{\lambda}\sum_{j}\int_{Q_{j}}|b(x)-b_{Q_{j}}||f_{Q_{j}}|\, M_{r}w_{j}(x)dx\\
 &\leq&\dfrac{c}{\lambda}\sum_{j}\inf_{Q_{j}}M_{r}w_{j}\int_{Q_{j}}|b(x)-b_{Q_{j}}||f_{Q_{j}}|\,dx\\
&\leq&\dfrac{c}{\lambda}\sum_{j}\dfrac{1}{|Q_{j}|}\int_{Q_{j}}|b(x)-b_{Q_{j}}|
\int_{Q_{j}}|f(y)|M_{r}w_{j}(y)\,dy\,dx\\
&\leq&c \|b\|_{BMO}\int_{\R^{n}}{{|f(x)|}\over{\lambda}}M_{r}w(x)\,dx.
\end{eqnarray*}
Observe that the constant $c$ is dimensional.

For $B_1$ we have by the generalized H\"older inequality (\ref{GHI}) and John-Nirenberg's theorem \eqref{John-Nirenberg-LLogL}
\begin{eqnarray*}
B_1&=&\dfrac{c}{\lambda}\sum_{j}\int_{Q_{j}}|b(x)-b_{Q_{j}}||f(x)|\,M_{r}w_{j}(x)dx\\
&\leq&\dfrac{c \|b\|_{BMO}}{\lambda}\sum_{j}\inf_{Q_{j}}M_{r}w_{j}
|Q_{j}|\|f\|_{L\log L,Q_{j}}.
\end{eqnarray*}
Now combining formula (\ref{RaoRen}) together with (\ref{CZdec}) and
recalling that $\Phi(t)= t (1+\log^{+}t)$ we have
\begin{eqnarray*}
\frac{1}{\lambda} |Q_{j}|\|f\|_{L\log L,Q_{j}}&\leq&
\frac{1}{\lambda} |Q_{j}|
\inf_{\mu>0}\{\mu+\dfrac{\mu}{|Q_{j}|}\int_{Q_{j}}\Phi\left(\dfrac{|f(x)|}{\mu}\right)dx\}\\
&\leq&|Q_{j}|+\int_{Q_{j}}\Phi\left(\dfrac{|f(x)|}{\lambda}\right)dx\\
&\leq& \frac{1}{\lambda} \int_{Q_{j}}|f(x)|\,dx +
\int_{Q_{j}}\Phi\left(\dfrac{|f(x)|}{\lambda}\right)dx\\ &\leq& 2
\int_{Q_{j}}\Phi\left(\dfrac{|f(x)|}{\lambda}\right)dx.
\end{eqnarray*}
Then
\begin{eqnarray*}
B_1&\leq &c\,\|b\|_{BMO}
\sum_{j}\int_{Q_{j}}\Phi\left(\dfrac{|f(x)|}{\lambda}\right)
\,w(x)dx\\ &\leq&
c\|b\|_{BMO}\int_{\R^{n}}\Phi\left(\dfrac{|f(x)|}{\lambda}\right)
\,w(x)dx.\\
\end{eqnarray*}

Combining this with estimates for $I$ and $II$, we get that
\begin{equation*}
w({x\in \R^n: |[b,T]f(x)|>\lambda})\leq
\dfrac{c\|b\|_{BMO}^{p}}{\lambda}
(pp')^{2p}[w]_{A_1}^{2p}\int_{\R^n}\Phi\left(\dfrac{|f(x)|}{\lambda}\right)\,M_{LlogL}w(x)dx.
\end{equation*}

Setting here $p=1+\frac{1}{\log(1+\|w\|_{A_1})}$ gives the estimate
we were looking for finishing the proof, because $ \frac{1}{\log(1+\|w\|_{A_1})}<1$ and the fact that for every $A>1$, $A^{1/A}<e$.

\hfill $\Box$

\subsection*{Acknowledgements}
I would like to thank my two advisors, Carlos P\'{e}rez and Rodrigo
Trujillo for our discussions, their support, help and friendship, and to the referee for the very useful remarks and corrections.

\end{document}